\documentstyle[din_a4,12pt]{article}
\begin{document}
\begin{center}
\Large
{\bf Generalized Steiner's Problem and its Solution with the Concepts in Field Thoery} \\
\large
\vspace{1.0cm}
Hai Lin \\
{Feb. 14th, 2001}
\\
 \vspace{0.5cm}
\normalsize
 $Department of Physics, Peking University, P.R.China, 100871$\\
{$ Email:hailin@mail.phy.pku.edu.cn$} \\
\vspace{0.5cm}
\large
{\bf Abstract} \\
\end{center}
\vspace{-0.3cm}
\normalsize

We generalized the Steiner's shortest line problem and found its connection 
with the concepts in classical field theory. We solved the generalized Steiner's problem by introducing
a conservative potential and a dissipative force in the field  and gave a computing method
by using a testing point and a corresponding iterative curve. 

\vspace{1.0cm}
\section{Steiner's Problem}
In the $ {\Re}_{}^{2} $ Euclidean space, there are $n$ points ${\alpha}_{i}^{} (i=1, ...,n.)$. 
Find a point ${\alpha}_{steiner}^{} $ that the sum of 
its distances to the $n$ points is a minimum: \\
\begin{equation}
\forall\alpha\in{\Re}_{}^{2} , {\sum}_{i=1}^{n} |{\alpha}_{steiner}^{}-{\alpha}_{i}^{}| 
\leq{\sum}_{i=1}^{n}  |\alpha-{\alpha}_{i}^{}| . 
\end{equation}
This is the original Steiner's problem[1].

\section{Generalized Steiner's Problem }
We would like to expand this problem to a more general situation: \\
In the D-dimension differentiable linear space ${E}_{}^{D} $, there are $n$ points ${\alpha}_{i}^{} (i=1, ...,n.)$. 
How to find a point ${\alpha}_{steiner}^{} $ that the sum of its defined distances $dis$ to the $n$ points is a minimum: \\
\begin{equation}
\forall\alpha\in{E}_{}^{D} , {\sum}_{i=1}^{n} dis ( {\alpha}_{steiner}^{}-{\alpha}_{i}^{}) 
\leq{\sum}_{i=1}^{n}  dis (\alpha-{\alpha}_{i}^{}) . 
\end{equation}
We call the point ${\alpha}_{steiner}^{} $  as the Steiner point.
 \section{Mapping Concepts onto  Field Theory}
If we map the concept of 'distance' in the generalized Steiner's problem onto the  concept of 'force potential'
in the Filed Theory, we will obtain another equivalent form of the generalized problem.
Let's consider that a particle at the testing point has conservative force interactions with each of the $n$
points and therefore has potentials  ${U}_{i}^{}(\alpha-{\alpha}_{i}^{})$ $ (i=1, ...,n.) $. 
The Steiner point makes the sum of the $n$ potentials to be a minimum.
If we let
\begin{equation}
{U}_{i}^{}(\alpha-{\alpha}_{i}^{})=dis(\alpha-{\alpha}_{i}^{})  ,  (i=1, ...,n.) 
\end{equation}
and find the Steiner point in this conservative field, the generalized Steiner's problem is solved. 
\\Obviously, the Steiner point should be in the set:
\begin{equation}
{\alpha}_{steiner}^{} \in set \left< \alpha | \nabla \left( {\sum}_{i=1}^{n} {U}_{i}^{}(\alpha-{\alpha}_{i}^{})\right)=0 \right>.
\end{equation}
Thus ${\alpha}_{steiner}^{} $ can be identified by comparing each elements in the set.

\section{Testing and Iterating Methods}
The next question is 'How to get that set by a computing method ?'  We would like to use a testing point $ {\alpha}_{testing}^{}$
and the interative curve ${L}_{iterative}^{}$ and introduce a conservative potential $U(\alpha)= {\sum}_{i=1}^{n} {U}_{i}^{}(\alpha-{\alpha}_{i}^{})
$ and a dissipative force. Note that the disspative force acts only when a particle moves.  
If $ {\alpha}_{testing}^{}$ is not in $ set \left< \alpha | \nabla \left( {\sum}_{i=1}^{n} {U}_{i}^{}(\alpha-{\alpha}_{i}^{})\right)=0 \right> $,
we put a point-like large-mass rest particle at  $ {\alpha}_{testing}^{}$. The conservative force exerted to the particle is[2]:
\begin{equation}
-\nabla U({\alpha}_{testing}^{}). 
\end{equation}
The particle is intended to move along the direction of the conservative force, but  the disspative force is strong enough to keep 
the particle moving in a quasi-static matter until rest at a point ${\alpha}_{s}^{}$. (${\alpha}_{s}^{} \in  set \left< \alpha | 
\nabla ({\sum}_{i=1}^{n} {U}_{i}^{}(\alpha-{\alpha}_{i}^{}))=0 \right>$.) \\
If the particle moves slowly enough, its path from ${\alpha}_{testing}^{}$ to ${\alpha}_{s}^{}$ is regarded as an iterative curve,
with which we can find ${\alpha}_{s}^{}$ by iterating points along the curve from ${\alpha}_{testing}^{}$.\\

\section{Iterative Curve and Conservative Potential}
The next question is 'How to get the iterative curve $L$ after the testing point is known?'. The iterative curve's tangent is everywhere parallel to the conservative force
 (For the D-dimemsion situation, the coordinates are ${Z}_{i}{}  ,(i=1,2,3,...D.)$ and $\lambda$ is a nonzero real number):

\begin{equation}
\left( 1,  {{\partial {Z}_{2}{}} \over {\partial {Z}_{1}{}}} {\huge |}_{L}{},{{\partial {Z}_{3}{}} \over {\partial {Z}_{1}{}}} {\huge |}_{L}{}, 
 ..., {{\partial {Z}_{D}{}} \over {\partial {Z}_{1}{}}} {\huge |}_{L}{} \right ) = {\lambda} \left( {{\partial U} \over {\partial {Z}_{1}{}}},
{{\partial U} \over {\partial {Z}_{2}{}}}, {{\partial U} \over {\partial {Z}_{3}{}}} ,..., {{\partial U} \over {\partial {Z}_{D}{}}} \right ).
\end{equation}
Thus we get:
\begin{equation}
 {{\partial {Z}_{i}{}} \over {\partial {Z}_{1}{}}} {\huge |}_{L}{} ={  {{\partial U} \over {\partial {Z}_{i}{}}} {\Huge /} 
{{\partial U} \over  {\partial {Z}_{1}^{}}}  }, ( i=2,3,...D ).
\end{equation}
If the testing point $ {\alpha}_{testing}^{}=({Z}_{1}^{testing}, {Z}_{2}^{testing}, ...,{Z}_{D}^{testing}) $,
the iterative curve is governed by the following D-1 functions:
\begin{equation}
 {Z}_{i}{}-{\int}_{{Z}_{1}^{testing}}^{{Z}_{1}^{}} \left(  { {{\partial U} \over {\partial {Z}_{i}{}}} {\Huge /} 
{{\partial U} \over  {\partial {Z}_{1}^{}}}  }  \right) d{Z}_{1}^{} =0  , ( i=2,3,...D ).
\end{equation}

\end{document}